\documentclass{article}
\usepackage{amsmath}
\usepackage{amsfonts}
\usepackage{amssymb}
\usepackage{amsthm}
\usepackage{graphicx}
\usepackage{cite}
\usepackage{url}
\usepackage{tikz}

\newtheorem{thm}{Theorem}[section]
\newtheorem{defi}{Definition}[section]
\newtheorem{remk}{Remark}[section]

\newtheorem{example}{Example}[section]
\newtheorem{propo}{Proposition}[section]

\newtheorem{notat}{Notation}[section]
\newtheorem{corol}{Corollary}[section]

\title{$ \infty $-CATEGORY AND SOME APPLICATIONS ON ORBIFORDS}
\author{JIAJUN DAI\thanks{jjdai@mail.ustc.edu.cn}}
\date{}
\begin{document}
	
	\maketitle
	
\section{Introduction}
	This paper in progress is motivated by a question about a given orbispace $ X $ whether there exists a global presentation $ X = Y / G $ for $ G $ a compact Lie group (such an orbispace is usually called a global quotient). \\
	L{\"u}ck–Oliver \cite{luck2001completion} made that if $ X = Y / \Gamma $ for a discrete group $ \Gamma $, then $ X $ is a global quotient; Henriques–Metzler \cite{henriques2004presentations} gave that every paracompact smooth effective $ n $-dimensional orbifold is a global quotient (of its orthonormal frame bundle by $ O(n) $), and a sufficient condition for a general (not necessarily effective) orbifold to be a global quotient; Henriques \cite{henriques2005orbispaces} conjectured that every compact orbispace is a global quotient; the analogous question for algebraic stacks has been studied by Edidin–Hassett–Kresch–Vistoli \cite{edidin2001brauer} and Totaro \cite{totaro2004resolution}; Kresch–Vistoli (\cite{kresch2004coverings} theorem2) (\cite{kresch2009geometry} theorem4.4) (using an unpublished result of Gabber) proved that smooth separated Deligne–Mumford stacks are global quotients; Pardon \cite{pardon2019enough} showed that all orbispaces satisfying very mild hypotheses are global quotients and particularly that all compact orbispaces are global quotients which verifies Henriques’ conjecture.\\
	Since an orbifold can be viewed as a groupoid which is also the fundamental groupoid of an $ \infty-category $, a naive idea to investigate orbifold lies in $ \infty-category $ theory \cite{kerodon} \cite{lurie2009higher}; on the other hand, as an orbifold is also a geometric/topological stack which takes values in groupoid, we could even extend higher stack and derived stack theories respectively to higher orbifold to obtain higher structure in orbifold and to derived orbifold to approximate an arbitrary but not too bad orbifold. Moreover, derived orbifold would become a stable $ \infty-category $, allowing us to study its stability via either stable $ \infty-category $ theory or stable homotopy theory or even equivariant homotopy theory.\\
	Here is mainly about an early result that orbifold stack is globally representable via some $ \infty-categorical $ techniques.
	
 \section{Preliminary $ \infty-category $ theory}
	In the setting of the category of topological spaces, homotopy could be intuitively viewed as a “path” between two continuous maps $ f_{1} $ and $ f_{2} $ between two topological spaces, recording that $ f_{1} $ could be “deformed” to $ f_{2} $ or the other way around, but not being able to tell us to what extent the “deformation” is. To obtain such information, a try could be made to regard the homotopy as a “path space”, and consider the “path”, which could also be viewed as homotopy and be called $ 2-homotopy $ to distinguish from the formal homotopy that we would call it $ 1-homotopy $ or $ 2-map $ (ordinary continuous map be called $ 1-map $), between two “path spaces” $ p_{1} $ and $ p_{2} $ which are literally homotopies, thus $ 2-homotopy $, called $ 3-map $ inductively, could form a new “2-path spaces”, encoding more information about the “deformation or derivation” between two $ 1-homotopies $ $ p_{1} $ and $ p_{2} $. Inductively, we could define $ n-homotopy $ which could be viewed as “n-path spaces” between two $ (n-1)-maps $ $ h_{1} $ and $ h_{2} $ to record the “derivation” between $ h_{1} $ and $ h_{2} $, where we call it the derivation of “degree $ n-1 $” in an extremely non-strict way. When $ n $ goes to $ \infty $, we might get all homotopies involved.\\
	As for approaches to realize the trial mentioned above, the Homotopy Hypothesis proposed by Grothendieck has been nominated  that the $ \infty-groupoids $ are equivalent to the topological spaces (considered modulo weak homotopy equivalence), and $ \infty-groupoids $ are equivalent to homotopy $ n-types $ for all extended natural numbers and moreover this equivalence is induced by the fundamental $ \infty-groupoid $ construction. A well-known model for $ \infty-groupoids $ in the context of higher category is $ Kan $ complex: specifically, $ Kan $ complex is homotopy equivalent to the singular complex of a topological space. Another approach to develop higher category is based on simplicial sets which serve as models for homotopy types: every $ CW $ complex is homotopy equivalent to the geometric realization of a simplicial set. Therefore, it is convenient to relate $ Kan $ complex with simplicial set to deal with homotopy theories combinatorically.
	
\begin{defi}
A \textbf{$ Kan $ complex} is a simplicial set $ K $ which has the following property: for any $ 0 \leq i \leq n $, any map $ f_{0} : \Lambda_{i}^{n} \rightarrow K $ admits an ($ Kan $) extension $ f : \Delta^{n} \rightarrow K $ or equivalently $ K_{n} \rightarrow \Lambda_{i}^{n} (K) $ is surjective for $ 0 \leq i \leq n $.
\end{defi}

\begin{example}
If $ X $ is a topological space, then its singular complex $ Sing(X) $ is a $ Kan $ complex; $ \mathcal{G} $ is a groupoid object, if and only if $ N(\mathcal{G}) $ is a $ Kan $ complex where $ N $ is the simplicial nerve functor.
\end{example}

For a category $ \mathcal{C} $, its simplicial nerve is characterized by the following formulation:\\
$ N(\mathcal{C})_{n} : \operatorname{Fun}( \overbrace{ \cdot \rightarrow \cdot \rightarrow ... \rightarrow \cdot}^{n+1 \: \text{copies}} , \mathcal{C}) = \left\lbrace X_{0} \rightarrow X_{1} \rightarrow ... \rightarrow X_{n} , X_{i} \in \mathcal{C} \right\rbrace $, whose $ i $th face map is obtained by composing the adjacent two morphisms in the ith position $ X_{i-1} \rightarrow X_{i} \rightarrow X_{i+1} $ and $ i $th degeneracy map by inserting the identity $ X_{i} \rightarrow X_{i} $. It is well known that there is a pair of adjoint functors $ | - | \dashv N $, where $ | - | $ is the geometric realization functor.

\begin{remk}
$ Kan $ complex is a simplicial model category which is enriched over the category of simplicial sets $ \text{Set}_{\Delta} $ endowed with Joyal model structure. To illustrate the model structure of $ Kan $ complex, the three distinguished classes of morphisms are:
\begin{itemize}
	\item cofibration: monomorphism which has the left lifting property with respect to all trivial $ Kan $ fibrations and whose collection is weak saturated (closed under pushout, retract and transfinite composition); acyclic/trivial cofibration: anodyne morphism which is defined to have the left lifting property with respect to all $ Kan $ fibrations.
	\item fibration: $ Kan $ fibration which has the right lifting property with respect to all horn inclusion $ \Lambda_{i}^{n} \rightarrow \Delta^{n} $ for all $ n > 0 $ , $ 0 \leq i \leq n $, and is closed under pullback, retract, composition and filtered colimit; acyclic/trivial fibration: trivial $ Kan $ fibration which has the right lifting property with respect to all boundary inclusion $ \partial \Delta^{n} \rightarrow \Delta^{n} $ for all $ n>0 $.
	\item weak equivalence: weak homotopy equivalence which could be easily checked to possess the ‘two out of three’ property.
\end{itemize}
\end{remk}

Thus, the weak factorization system of $ Kan $ complex is (anodyne morphisms, $ Kan $ fibrations) and (monomorphisms, trivial $ Kan $ fibrations).\\
This model structure implies that the fibrant objects of $ \text{Set}_{\Delta} $ are precisely the $ \infty-categories $.

\begin{defi}
An \textbf{$ \infty-category $} is a simplicial set $ K $ which has the following property: for any $ 0 < i < n $, any map $ f_{0} : \Lambda^{n}_{i} \rightarrow K $ admits an (weak $ Kan $) extension $ f : \Delta^{n} \rightarrow K $ or equivalently $ K_{n} \rightarrow \Lambda^{n}_{i}(K) $ is surjective for $ 0 < i < n $.
\end{defi}

Particularly, objects of an $ \infty-category $ are the $ 0-simplices $; the morphisms are the $ 1-simplices $; the horn-fillers exhibit not only the homotopies, i.e. horn-fillers of higher degree represent higher homotopies, but also the associator and unit of the composition law in higher category. Furthermore, for any objects $ x,y \in \mathcal{C} $, the mapping space $ \operatorname{Map}_{\mathcal{C}}(x,y) $ is the $ Kan $ complex whose $ n-simplices $ are maps $ \Delta^{n} \times \Delta^{1} $ to $ \mathcal{C} $ which sends $ \Delta^{n} \times \left\lbrace 0 \right\rbrace $ to the vertex $ x $ and $ \Delta^{n} \times \left\lbrace 1 \right\rbrace $ to the vertex $ y $.

\begin{remk}
$ \infty-category $, as named by Lurie, is literally, or more precisely a model for, $ (\infty,1)-category $ (whose $ k-morphisms $ are invertible for $ k \geq 2 $) which was called quasi-category by Joyal and first referred to as weak $ Kan $ complex by Boardman-Vogt. Other models for $ (\infty,1)-category $ like Segal category, Rezk category (also known as complete Segal category) are equivalent to $ \infty-category $. These models are employed to keep track of all homotopies as demonstrated before.
\end{remk}

\begin{example}
Obviously, any $ Kan $ complex is an $ \infty-category $, so is $ Sing(X) $ for a topological space $ X $; for an ordinary category $ \mathcal{C} $, $ N(\mathcal{C}) $ is an $ \infty-category $, moreover, $ S $ is isomorphic to the nerve of an ordinary category $ \mathcal{C} $ if and only if, for each $ 0 < i < n $, the map $ S_{n} \rightarrow \Lambda^{n}_{i}(S) $ bijective; the Eilenberg-Mac Lane space $ K(A, n) $ for an abelian group $ A $ is an $ \infty-category $ for $ n \geq 1 $.
\end{example}

\begin{propo}[\cite{kerodon}]
For an $ \infty-category $ (respectively $ Kan $ complex) $ \mathcal{C} $, for a simplicial set $  _{\cdot} $, $ \operatorname{Fun}(S_{\cdot}, C) $ is again an $ \infty-category $(respectively $ Kan $ complex).
\end{propo}
	
In particular, if $ \mathcal{D} $ is an $ \infty-category $ (respectively $ Kan $ complex), $ \operatorname{Fun}(\mathcal{D}, \mathcal{C}) : = \operatorname{Map} _{\text{Set}_{\Delta}} (\mathcal{D}, \mathcal{C}) $ is an $ \infty-category $ (respectively $ Kan $ complex), which is referred to as the $ \infty-category $ of functors and whose set of $ n-simplices $ are by definition $ \operatorname{Hom} _{\text{Set}_{\Delta}} (\mathcal{D} \times \Delta^{n}, \mathcal{C}) $. This property ensures that $ \infty-category $ and $ Kan $ complex are internal. 
	
\begin{defi}
For an $ \infty-category \: \mathcal{C} $, an object $ X \in \mathcal{C} $ is \textbf{ final }(respectively \textbf{initial}) if for any $ Y \in \mathcal{C} $, $ \operatorname {Map}_{\mathcal{C}}(Y,X) $ (respectively $ \operatorname {Map}_{\mathcal{C}}(Y,X) $) is contractible.
\end{defi}

\begin{defi}
For simplicial sets $ S, S' $, their \textbf{join} $ S \star S' $ is characterized by the formulation $ (S \star S')_{n} : S_{n} \bigcup S'_{n} \bigcup \cup_{i+j=n-1} S_{i} \times S'_{j} $.
\end{defi}

It can be easily seen that $ \Delta^{i} \star \Delta^{j} $ is isomorphic to $ \Delta^{i+j+1} $. Moreover, for two ordinary categories $ \mathcal{C} $ and $ \mathcal{C'} $, $ N(\mathcal{C} \star \mathcal{C'}) $ is isomorphic to $ N(\mathcal{C}) \star N(\mathcal{C'}) $.

\begin{propo}[Joyal]
The join of two $ \infty-categories $ is an $ \infty-category $.
\end{propo}

\begin{notat}
For simplicial set $ K $, the left cone $ K^{\vartriangleleft} $ is defined to be the join $ \Delta^{0} \star K $ . Dually, the right cone $ K^{\vartriangleright} $ is defined to be the join $ K \star \Delta^{0} $. Either cone contains a distinguished vertex (belonging to $ \Delta^{0} $ ), which will be referred to as the cone point.
\end{notat}

\begin{propo}[Joyal]
If $ p : K \rightarrow S $ is a map of simplicial sets, then there exists a simplicial set $ S _{/p} $ with the following universal property: $ \operatorname{Hom}_{\text{Set}_{\Delta}}\left(Y, S_{/ p}\right)=\operatorname{Hom}_{p}(Y \star K, S) $, where the subscript on the right hand side means that the hom subset consists of only morphisms $ f: Y \star K \rightarrow S $ such that $ f|_{K} = p $.
\end{propo}

\begin{notat}
For above map $ p $, if S is an $ \infty-category $, \textbf{$ S_{/ p} $} denotates an overcategory of $ S $ or the $ \infty-category $ of objects of $ S $ over $ p $; dually, replacing $ Y \star K $ by $ K \star Y $, \textbf{$ S_{p /} $} denotes an undercategory of $ S $ or the $ \infty-category $ of objects of $ S $ under $ p $.
\end{notat}

\begin{defi}[Limits and colomits]
Let $ \mathcal{C} $ be an $ \infty-category $ and let $ p: K \rightarrow \mathcal{C} $ be an arbitrary map of simplicial sets. A colimit for $ p $ is an initial object of $ \mathcal{C}_{p/} $ , and a limit for p is a final object of $ \mathcal{C}_{/p} $ .
\end{defi}

A colimit diagram is a map $ \bar{p} : K^{\vartriangleright} \rightarrow \mathcal{C} $ which is a colimit of $ P = \bar{p} |_{K} $ . In this case,  $ \bar{p} (\infty) \in \mathcal{C} $ is referred to as a colimit of $ p $ where $ \infty $ denotes the cone point of $ K^{\vartriangleright} $; similarly, a limit diagram is a map $ \bar{p} : K^{\vartriangleleft} \rightarrow \mathcal{C} $ which is a limit of $ P = \bar{p} |_{K} $, notated $ \bar{p} (-\infty) \in \mathcal{C} $  where $ - \infty $ denotes the cocone point of $ K^{\vartriangleleft} $.

\section{Orbifolds}

An orbifold is locally the quotient of a manifold by a finite group. A more general notion is the orbispace which is locally the quotient of a space by a group. There are two conventional approaches to orbifolds: one is via (topological) groupoid objects, and the other is via (topological) stacks. Both orbifold groupoid and orbifold stack can be added new ingredients of and reinterpreted by $ \infty - category $ theory.\\
We assume that every orbifold $ \mathcal{X} $ of concern has finite stablizers in $ \mathcal{F} $ and fix a family $ \mathcal{F} $ of allowed isotropy groups that is an essentially small class of paracompact topological groups which are closed under isomorphisms.

\begin{defi}
A (left) \textbf{action of groupoid} $ \mathbb{G} $ on a space $  X \rightarrow \mathbb{G}_{0} $ over $ \mathbb{G}_{0} $ is exhibited by a map a:  $ \mathbb{G}_{1} \times_{\mathbb{G}_{0}} X \rightarrow X $ satisfying the associative and unital conditions that the following diagrams commute
\[
\begin{array}{ccc}
	\mathbb{G}_{1} \times_{\mathbb{G}_{0}} \mathbb{G}_{1} \times_{\mathbb{G}_{0}} X & \stackrel{m \times 1}{\longrightarrow} & \mathbb{G}_{1} \times_{\mathbb{G}_{0}} X \\
	\downarrow^{1 \times a} & & \downarrow^{a} \\
	\mathbb{G}_{1} \times_{\mathbb{G}_{0}} X & \stackrel{a}{\longrightarrow} &  X
\end{array} 
\]
\[
\begin{array}{ccc}
	\mathbb{G}_{0} \times_{\mathbb{G}_{0}} X & \stackrel{i \times 1}{\longrightarrow} & \mathbb{G}_{1} \times_{\mathbb{G}_{0}} X \\
	\downarrow^{\simeq} & & \downarrow^{a} \\
	X & \stackrel{=}{\longrightarrow} &  X
\end{array}
\]
\end{defi}

We can briefly define this action groupoid as $ \mathbb{G} \ltimes X : = (\mathbb{G}_{1} \times _{\mathbb{G}_{0}} X \rightrightarrows X) $. Such space $ X \rightarrow \mathbb{G}_{0} $ (or briefly referred to as $ X $) equipped with a (left) action of $ \mathbb{G} $ is called a $ \mathbb{G} - space $.\\
Roughly speaking, spaces over $ \mathbb{G}_{0} $ could be viewed as a continuous $ \mathbb{G}_{0} - indexed $ family of spaces and a left action of $ \mathbb{G} $ on a space over $ \mathbb{G}_{0} $ as a continuous contravariant functor from $ \mathbb{G} $ to spaces.\\
For a space $ U \rightarrow  \mathbb{G}_{0} $, the \textbf{restriction} $ \mathbb{G}_{U} $ of $ \mathbb{G} $ to $ U $ consists of object space $ U $ and arrow space $ U \times _{\mathbb{G}_{0}} \mathbb{G}_{1} \times _{\mathbb{G}_{0}} U $.

\begin{defi}
A \textbf{cellular groupoid} $ \mathbb{G} $ is a topological groupoid of the form $ \text{colim} \: \mathbb{G}_{\alpha} $ whose structure is built as follows:
\[ 
\begin{array}{ccc}
	S^{n-1} \times \mathbb{O} & \longrightarrow & D^{n} \times \mathbb{O}\\
	\downarrow^{f} & & \downarrow \\
	\mathbb{G}_{\alpha -1} & \longrightarrow &  \mathbb{G}_{\alpha}
\end{array}
 \]
where $ \mathbb{G}_{n} : = \operatorname{map}(\Delta^{n}, \mathbb{G}) $; $ \alpha $ runs through the set of ordinals smaller than a given ordinal; $ \mathbb{O} $ is an orbit groupoid with the action groupoid of the form $ G \ltimes (G/H) $ for $ G/H $ being a $ G-orbit $ with allowed isotropy group; $ \mathbb{G}_{\alpha} $ is obtained from $ \mathbb{G}_{\alpha-1} $ via the pushout of the cell groupoid $ D^{n} \times \mathbb{O} $ along an attaching map $ f $ if $ \alpha-1 $ exists, and as the colimit $ \varinjlim_{\beta < \alpha} \mathbb{G}_{\beta} $ otherwise. 
\end{defi}

\begin{defi}
Given a topological groupoid $ \mathbb{G} $ , a subspace $ Z \subset \mathbb{G}_{0} $ is called \textbf{saturated} if there are no arrows from $ Z $ to its complement. A \textbf{saturated inclusion} of topological groupoids is a map of the form $ \mathbb{G}_{Z} \rightarrow \mathbb{G} $ such that $ Z \subset \mathbb{G}_{0} $ is a saturated subspace and $ \mathbb{G}_{Z} $ is the restriction of $ \mathbb{G} $ to $ Z $.
\end{defi}

The map $ S^{n-1} \times \mathbb{O} \rightarrow D^{n} \times \mathbb{O} $ associated with a cellular groupoid is a closed saturated inclusion which can easily be checked to be compatible with cofibrations in simplicial model category $ \text{Set}_{\Delta} $.\\
Define a cofibrant groupoid to be a retract of a cellular groupoid, it is still cellular.

\begin{propo}[Gepner-Henriques\cite{gepner2007homotopy}]
Given topological $ \mathbb{H} $ and $ \mathbb{G} $ , there is a natural equivalence of groupoids
\[ 
\{ Principal \mathbb{G}-bundles on \mathbb{H} \} \simeq \operatorname{hocolim}_{U \in \operatorname{Cov}(\mathbb{H}_{0})} \operatorname{\mathbb{H}om}(\mathbb{H}_{U}, \mathbb{G})
 \],
where $ \operatorname{Cov}(\mathbb{H}_{0}) $ denotes the cover of $ \mathbb{H}_{0} $.
\end{propo}

This proposition tells us that we could recover the principal bundle on a certain groupoid $ \mathbb{H} $ that has a global section by gluing local sections of this principal bundle with respect to a well-refined(cofinal) cover of the underlying space $ \mathbb{H}_{0} $ of $ \mathbb{H} $.

\begin{defi}
A topological groupoid $ \mathbb{G} $ is \textbf{fibrant} if $ \mathbb{G}_{1} \rightarrow \mathbb{G}_{0} $ is a universal $ \mathbb{G}-bundle $. A retract of a fibrant groupoid is still fibrant.
\end{defi}

\begin{defi}
The \textbf{fibrant replacement} $ \text{fib}(\mathbb{G}) $ of a topological groupoid $ \mathbb{G} $ is the gauge groupoid of the universal principal $ \mathbb{G} -bundle \: || \mathbb{EG} || \rightarrow || \mathbb{G} || $; that is, $ \text{fib}(\mathbb{G}):=\text{gauge }(||\mathbb{EG}||) = ((||\mathbb{EG}|| \times _{\mathbb{G}_{0}}||\mathbb{EG}||)/\mathbb{G} \rightrightarrows ||\mathbb{G}|| )$ , where $ ||-|| $ is the fat geometric realization functor, $ \mathbb{EG} $ the translation groupoid associated to $ \mathbb{G} $ with the action groupoid $ \mathbb{G} \ltimes \mathbb{G}_{1} $. 
\end{defi}

Fibrant replacement preserves pushouts along closed saturated inclusions which means it behaves well in the operation of gluing cells. Specifically, the fibrant replacement of a cellular groupoid is cellular.

\begin{propo}[Gepner-Henriques\cite{gepner2007homotopy}]
Let $ \mathbb{G} $ be any topological groupoid. Then $ \text{fib}(\mathbb{G}) $ is fibrant.
\end{propo}

Furthermore, we can adapt cellular groupoid to be simplicial groupoid.

\begin{defi}
Let $ \mathbb{G}^{n} : = \overbrace{ \mathbb{G}_{1} \times _{\mathbb{G}_{0}} \mathbb{G}_{1} \times _{\mathbb{G}_{0}} ... _{\mathbb{G}_{0}} \mathbb{G}_{1}}^{n \: \text{copies}} , \mathcal{C}) $ , then $ \mathbb{G}^{n} \hookrightarrow \mathbb{G}^{n+1} $ are saturated inclusions. In this case, we can define $ \mathbb{G} : = \text{colim} (\mathbb{G}^{1} \hookrightarrow \mathbb{G}^{2} \hookrightarrow ...) $ to be a \textbf{simplicial groupoid}.
\end{defi}

A simplicial groupoid is naturally cellular, thus its cofibrant object and fibrant replacement could be defined in the same fashion as those of a cellular groupoid and inherit the latters’ properties.\\
On the other hand, since a simplicial groupoid is by defibition an $ \infty -category $, it is fibrant subsequently.

\begin{defi}
A lax presheaf of groupoids $ \mathcal{X} := \operatorname{Fun}(Top^{op}, Grpd) $ on $ Top $, the category of topological spaces, is a \textbf{topological stack} if, for each object $ T $ of Top and each cover $ U \rightarrow T $ of $ T $, the natural map
\[ 
\mathcal{X}(T) \rightarrow \operatorname{holim} \left\{\mathcal{X}(U) \rightrightarrows \mathcal{X}\left(U \times_{T} U\right) 
\begin{tikzpicture}
\draw [->] (0,0)--(0.3,0);
\draw [->] (0,0.1)--(0.3,0.1);
\draw [->] (0,0.2)--(0.3,0.2);
\end{tikzpicture}
\mathcal{X}\left(U \times_{T} U \times_{T} U\right)\right\} 
 \]
is an equivalence of groupoids.
\end{defi}

Let a functor $ \mathcal{M} : \{ \text{Topological groupoids} \} \rightarrow \{ \text{Topological stacks} \} $ send the topological groupoid $ \mathbb{G} $ to the stack $ \mathcal{M}_{\mathbb{G}} : \{ \text{Topological spaces} \} \rightarrow \{ \text{Groupoids} \} $ which is defined by $ \mathcal{M}_{\mathbb{G}}(T) : = \{ \text{Principal } \mathbb{G} \text{ -bundles on } T \} $.\\

It can be easily seen that The functor $ \mathcal{M}_{\mathbb{G}} $ is the stackification of the presheaf of groupoids $ \operatorname{Hom} (-, \mathbb{G} ) $ and that the functor $ \operatorname{Hom} (-, \mathbb{G} ) $ classifies trivial principal $ \mathbb{G} -bundles $ up to isomorphism. What’s more, if $ \mathbb{G} $ is fibrant, then $ \operatorname{Hom} (-, \mathbb{G} ) $ is a stack.

\begin{propo}[Gepner-Henriques\cite{gepner2007homotopy}]
Let $ \mathbb{H} $ be a topological groupoid with paracompact object space. Then for any topological groupoid $ \mathbb{G} $, there is a natural equivalence of groupoids $ \operatorname{Hom} (\mathbb{H}, \text{fib} \mathbb{G}) \simeq \mathcal{M}_{\mathbb{G}} (\mathbb{H}) $. 
\end{propo}

Roughly speaking, we can simplify this proposition with a slogan that fibrant topological groupoids are stacks. For simplicial groupoids are fibrant, we can soon get the result:

\begin{corol}
Let $ \mathbb{H}, \mathbb{G} $ be a pair of simplicial groupoids where $ \mathbb{H} $ is paracompact, $ \operatorname{Hom} (\mathbb{H}, \text{fib} \mathbb{G}) \simeq \mathcal{M}_{\mathbb{G}} (\mathbb{H}) $.  
\end{corol}

\begin{thm}
For $ \mathbb{G} $ an orbifold groupoid, $ \mathcal{M}_{\mathbb{G}} $ is globally representable.
\end{thm}

Sketch of the proof\\
Step 1: Slice\\
For a stabilizer $ H_{x_{i}} $ in $ FlowerF $, where $ x_{i} $ is the corresponding fixed point, there is an action groupoid $ \mathbb{G} \ltimes H_{x_{i}} $, denoted $ _{i}\mathbb{H}  $. Then $ _{i}\mathbb{H} $ can be constructed by gluing cells one by one in finite times for $ H_{x_{i}} $ is finite thus can be obtained in a simplicial way, i.e. it is a simplicial groupoid. By corollary 3.1, for each action groupoid $ _{i}\mathbb{H}  $ , there is a (sub)stack $ \mathcal{M}_{\mathbb{G}}(_{i}\mathbb{H}) $.\\
Step 2: Join\\
By definition, $ \mathcal{M}_{\mathbb{G}}(_{i}\mathbb{H}) $is principal $ \mathbb{G} -bundles $ on $ (_{i} \mathbb{G} )_{0} $ with an (left) action $ \mathbb{G} \ltimes Hx_{i} $, thus can be recovered with respect to a collection of cofinal covers of $ (_{i} \mathbb{H})_{0} $ by proposition 3.1. Since $ _{i}\mathbb{H}  $ is simplicial, we can reason in two directions: simplicial groupoids $ _{i}\mathbb{H}  $ for all $ i $ can be joined to obtain a new simplicial groupoid $ \mathbb{H} $; such covers $ \operatorname{Cov}((_{i} \mathbb{H})_{o}) $ need to be simplicial thus cofinal covers become $ \infty -categories $. Join $ \operatorname{Cov}((_{i} \mathbb{H})_{o}) $ for all $ i $ whose number is finite to get a cover $ \operatorname{Cov} $ of the underlying object space $ \mathbb{H}_{0} $, then $ \operatorname{Cov} $ would be again $ \infty-category $ by proposition 2 to be cofinal. Hence, we can join substacks $ \mathcal{M}_{\mathbb{G}} (_{i}\mathbb{H}) $ to obtain a new stack $ \mathcal{M}_{\mathbb{G}}(\mathbb{H}) $ by proposition 3.1.

\bibliographystyle{alpha}
\bibliography{ref}

\end{document}